\documentclass[fleqn,10pt]{IEEEtran}
\pdfoutput=1

\usepackage{graphicx}
\usepackage{amssymb}
\usepackage{amsmath}
\usepackage{amsthm}
\usepackage{epstopdf}
\usepackage{enumerate}

\DeclareMathOperator{\std}{std}

\newcommand{\R}{\ensuremath{\mathbb{R}}}

\theoremstyle{definition}
\newtheorem{defn}{Definition}
\newtheorem{lem}{Lemma}
\newtheorem{prop}{Proposition}
\newtheorem{experiment}{Experiment}

\title{The Douglas-Rachford Algorithm for Weakly Convex Penalties}
\author{\.Ilker Bayram and Ivan W. Selesnick 
\thanks{\.{I}. Bayram is with the Dept. of Electronics and Telecommunications Engineering, Istanbul Technical University, Istanbul, Turkey (e-mail~:~ibayram@itu.edu.tr).} \thanks{ I.~W.~Selesnick is with the Dept. of Electrical and Computer Engineering, Polytechnic School of Engineering, New York University, NY, USA (e-mail~:~selesi@poly.edu).}}
\date{}

\begin{document}
\maketitle
\begin{abstract}
The Douglas-Rachford algorithm is widely used in sparse signal processing for minimizing a sum of two convex functions. In this paper, we consider the case where  one of the functions is weakly convex but the other is strongly convex so that the sum is convex. We provide a condition that ensures the convergence of the same Douglas-Rachford iterations, provided that the strongly convex function is smooth. We also present a modified Douglas-Rachford algorithm that does not impose a smoothness condition for the convex function. We then provide a discussion on the convergence speed of the  two types of algorithms and demonstrate the discussion with numerical experiments.
\end{abstract}
\section{Introduction}
The Douglas-Rachford algorithm is a widely used splitting method for solving problems of the form
\begin{equation}\label{eqn:problem}
\min_x \bigl\{ h(x) = f(x) + g(x) \bigr\}, 
\end{equation}
where $f$\, and $g$ are convex functions \cite{combettes_chp}. The algorithm employs proximity operators of $f$ and $g$ (see Defn.~\ref{def:prox} below), which are typically easier to realize than the proximity operator of the sum $h = f + g$. In this paper, we study the case where $g$ is weakly convex and $f$ is strongly convex so that their sum $h$ is convex. 
First, we show that the regular Douglas-Rachford iterations, that treat $f$ and $g$ as if they are convex functions, converge provided a condition on the step-size is satisfied. Second, we derive another algorithm by adding and subtracting a quadratic from $g$ and $f$ respectively so that the resulting functions are convex. In order to compare the algorithms, we study their convergence speed for a special case and demonstrate our findings via numerical experiments.

\subsection*{Douglas-Rachford Algorithm for a Convex Pair of Functions}
Before we discuss  weakly-convex penalties, let us consider the Douglas-Rachford algorithm for convex functions.
The algorithm can be succintly described in terms of the proximity operators associated with $f$ and $g$.
\begin{defn}\label{def:prox}
The proximity operator $J_{\alpha f}$ of a function ${f: \R^n \to \R}$ is defined as
\begin{equation}\label{eqn:defprox}
J_{\alpha f }(x) = \arg \min_{z} \frac{1}{2\alpha} \|z -x \|_2^2 + f(z).
\end{equation}
\end{defn}

If $\|x\|_2^2/(2 \alpha) + f(x)$\, is strictly convex, then the minimum of the problem in \eqref{eqn:defprox} is unique and $J_{\alpha f}$\, is well-defined.

Given this definition, the Douglas-Rachford iterations \cite{eck92p293,combettes_chp} for the problem in \eqref{eqn:problem} are
\begin{equation}\label{eqn:DR}
z^{n+1} = \Bigl (1-\lambda)\, I + \lambda\,\bigl(2J_{\alpha f} - I \bigr) \, \bigl(2J_{\alpha g} - I \bigr) \Bigr)\,(z^n).
\end{equation}
When both $f$ and $g$ are convex, these iterations produce a convergent sequence where the limit is related to a minimizer of $h$. More precisely,
\begin{prop}\cite{lio79p964, eck92p293}\label{prop:convex}
Suppose $f$ and $g$ are convex, and the set of minimizers of $h = f+g$ is non-empty.  For ${\alpha > 0}$ and $0<\lambda<1$, the sequence constructed by \eqref{eqn:DR} converges to a point $z^*$ such that $J_{\alpha g}(z^*)$ minimizes $h$.
\end{prop} 

We remark that $f$ and $g$ are assumed to be merely convex. Therefore, the roles of $f$ and $g$ are interchangeable in \eqref{eqn:DR} and Prop.~\ref{prop:convex}.

\subsection*{Douglas-Rachford Algorithm in a Weakly Convex Setting}
In sparse signal processing, one of the functions (say $g$) usually plays the role of a sparsity inducing penalty by means of its proximity operator (also referred to as a shrinkage/threshold function). However, when $g$ is convex, this threshold function introduces bias in the non-zero estimates (see e.g. Thm.~1 in \cite{sel14p078}). In order to circumvent such bias, non-convex functions are viable alternatives because they penalize high-magnitude coefficients less compared to convex penalties, which induces sparsity more effectively  \cite{sel14p078,bay15p265,cha14ICASSP}. In this paper, we restrict our attention to a special family of non-convex penalty functions named weakly-convex. 
\begin{defn}\cite{via83p231}
$g:\R^n \to \R$  is said to be $\rho$-weakly convex if 
\begin{equation}
g(x) + \frac{s}{2}\,\|x\|_2^2
\end{equation}
is convex for $s\geq \rho$.
\end{defn}
Weakly-convex penalties are of interest because the proximity operator of this family includes any separable monotone threshold function \cite{bay15p265,cha14ICASSP,ant07p16}.

When $g$ is not convex but $\rho$-weakly convex, convergence of the iterations in \eqref{eqn:DR} depend on the value of $\alpha$. Specifically, we show in this paper the following result.
\begin{prop}\label{prop:main}
Suppose $g$ is $\rho$-weakly convex, $f(x) - \dfrac{\rho}{2} \|x\|_2^2$ is convex, and $\nabla f$ is Frechet differentiable with Lipschitz constant $\sigma$, i.e., 
\begin{equation}
\| \nabla f(x) - \nabla f(y) \|_2 \leq \sigma \|x - y \|_2.
\end{equation}
Assume also that the set of minimizers of $h = f+g$ is non-empty.
If $\alpha \leq 1/ \sqrt{\sigma\,\rho}$ and $0< \lambda< 1$, then,
\begin{enumerate}[(a)]
\item \label{alg:a} the sequence produced by
\begin{equation}\label{eqn:DR1a}
z^{n+1} = \Bigl( (1-\lambda)\,I + \lambda\,\bigl(2J_{\alpha f} - I \bigr) \, \bigl(2J_{\alpha g} - I \bigr) \Bigr)\,(z^n),
\end{equation}
converges to some $z^*$ such that $J_{\alpha g}(z^*)$ minimizes $h$;
\item \label{alg:b} the sequence produced by
\begin{equation}\label{eqn:DR1b}
z^{n+1} = \Bigl( (1-\lambda)\,I + \lambda\,\bigl(2J_{\alpha g} - I \bigr) \, \bigl(2J_{\alpha f} - I \bigr) \Bigr)\,(z^n),
\end{equation}
converges to some $z^*$ such that $J_{\alpha f}(z^*)$ minimizes $h$. \qed
\end{enumerate} 
\end{prop}

We remark that, because $g$ is non-convex, there are fundamental differences between the properties of $J_{\alpha f}$ and $J_{\alpha g}$. Therefore, in contrast to the convex case, the two algorithms in Prop.~\eqref{prop:main} deserve to be studied separately.

One desirable feature of Prop.~\ref{prop:main} is that it generalizes Prop.~\ref{prop:convex} in the following sense. As $\rho \to 0$, i.e., as $g$\, approaches a convex function, the upper bound restricting $\alpha$ increases without limit. This is consistent with the convergence result in the convex case where there is no restriction on $\alpha$.

Prop.~\ref{prop:main} requires that $f$ be smooth on all $\R^n$, ruling out the inclusion of characteristic functions in $f$. A proposition with less stringent conditions on $f$, which also involves the proximity operators of $f$ and $g$ is given in the following.
\begin{prop}\label{prop:shift}
Suppose $g$ is $\rho$-weakly convex, $f(x) - \dfrac{\rho}{2} \|x\|_2^2$ is convex, the set of minimizers of $h = f+g$ is non-empty and $\alpha$ is a constant such that ${\alpha < 1/\rho}$. Also, let
\begin{equation}
\beta_1 = \frac{\alpha}{ 1 + \alpha\,\rho}, \qquad \beta_2 = \frac{\alpha}{ 1 - \alpha\,\rho},
\end{equation}
and 
\begin{equation}\label{eqn:K1K2}
K_1(x) = J_{\beta_1 g}\,\left( \frac{\beta_1}{\alpha} x \right), \qquad K_2(x) = J_{\beta_2 f}\,\left( \frac{\beta_2}{\alpha} x \right).
\end{equation}
Finally, suppose $0 < \lambda < 1$.
\begin{enumerate}[(a)]
\item If
\begin{equation}\label{eqn:DRKa}
z^{n+1} = \Bigl( (1-\lambda) I + \lambda\,\bigl(2K_2 - I \bigr) \, \bigl(2 K_1 - I \bigr) \Bigr)\,(z^n),
\end{equation}
then, the sequence of $z^n$'s  converge to some $z^*$ such that $K_1(z^*)$\, minimizes $h$. 
\item If
\begin{equation}\label{eqn:DRKb}
z^{n+1} = \Bigl( (1-\lambda) I + \lambda\,\bigl(2K_1 - I \bigr) \, \bigl(2 K_2 - I \bigr) \Bigr)\,(z^n),
\end{equation}
then, the sequence of $z^n$'s  converge to some $z^*$ such that $K_2(z^*)$\, minimizes $h$. \qed
\end{enumerate} 
\end{prop}

In contrast to Prop.~\ref{prop:main}, Prop.~\ref{prop:shift} allows $f$ to be non-smooth -- in particular $f$ may contain characteristic functions. Therefore, Prop.~\ref{prop:shift} applies to a larger class of splittings. 
Even though this remark seems to favor the algorithms in Prop.~\ref{prop:shift}, 
we will show in Section~\ref{sec:speed} that the algorithms in Prop.~\ref{prop:main} can be faster than the algorithms in Prop.~\ref{prop:shift}, under certain conditions.

\subsection*{Related Work and Contribution}
Proof of convergence for the Douglas-Rachford algorithm when both $f$ and $g$ are convex can be found in \cite{Bauschke,eck92p293,lio79p964}. The convergence proof in these works rely on a study of the mapping that takes $x^k$ to $x^{k+1}$ in \eqref{eqn:DR} and is not directly related with a monotone decrease of the cost at each iteration. In recent work, the authors of \cite{pat14CDC} construct an envelope function which monotonically decreases at each iteration, allowing one to interpret the Douglas-Rachford iterations as a variable metric gradient algorithm. A merit function based on this envelope function is used in \cite{li15arXiv} to study the Douglas-Rachford algorithm where 
$g$ is closed and $f$ is weakly convex. However, since these are very mild assumptions, general convergence can only be claimed provided that the Douglas-Rachford iterations produce a bounded sequence. Nevertheless, in \cite{li15arXiv}, the authors  show  that the algorithm converges for specific cases. The convergence condition in \cite{li15arXiv}  (see Thm.~1 of \cite{li15arXiv}), requires that $\alpha < (\sqrt{3/2} - 1) / \sigma$. 

Outside the convex setting, the convergence of the Douglas-Rachford algorithm has also been studied in \cite{hes13p397,borwein_chp, ara13p753} for non-convex feasibility problems. However, these works do not directly cover minimization problems considered in this paper.

In contrast to \cite{li15arXiv}, we work in a more restricted setting (strongly convex $f$, weakly convex $g$) and consequently obtain a stronger convergence statement (see Prop.~\ref{prop:main}). Specifically, Prop.~\ref{prop:main} ensures that the algorithm converges if  $\alpha < 1/\sqrt{\sigma \rho}$, allowing higher values of $\alpha$ than those in \cite{li15arXiv}. Our approach does not rely on the merit or envelope functions proposed in \cite{li15arXiv, pat14CDC} but on a study of the mapping that takes $x^k$ to $x^{k+1}$ as in \cite{Bauschke,eck92p293,lio79p964}. In the weakly convex case, a study of convergence following a similar plan was presented in \cite{bay15arXiv} for the iterative shrinkage/thresholding algorithm (ISTA).
However, the discussion in \cite{bay15arXiv} cannot be straightforwardly extended to cover the Douglas-Rachford algorithm and a study of the Douglas-Rachford algorithm in this setting has not appeared in the literature as far as we are aware.

\subsection*{Outline}
In the following, we present proofs of Prop.~\ref{prop:main} and Prop.~\ref{prop:shift} in Sections~\ref{sec:main}, \ref{sec:shift} respectively. For a special case, we study the convergence speed of the algorithms in Section~\ref{sec:speed} and demonstrate this discussion via experiments in Section~\ref{sec:exp}. Finally, we provide an outlook in Section~\ref{sec:conc}.

\section{Douglas-Rachford Iterations for Weakly Convex Penalties (Proof of Prop.~\ref{prop:main})}\label{sec:main}

Our  convergence analysis 
for Prop.~\ref{prop:main} relies on a study of the operators
\begin{subequations}\label{eqn:Ufgf}
\begin{align}
U_{f,g} &= \bigl(2J_{\alpha f} - I \bigr) \, \bigl(2J_{\alpha g} - I \bigr), \\
U_{g,f} &= \bigl(2J_{\alpha g} - I \bigr) \, \bigl(2J_{\alpha f} - I \bigr).
\end{align}
\end{subequations}
We show that 
\begin{enumerate}[(i)]
\item $U_{f,g}$ and $U_{g,f}$ become non-expansive under the conditions stated in Prop.~\ref{prop:main}. 
\item Fixed points of $U_{f,g}$ and $U_{g,f}$ are such that the image of these sets under $J_{\alpha f}$ and $J_{\alpha g}$ respectively, gives the set of minimizers of $f + g$. 
\end{enumerate}
Once these are shown, convergence follows by the Krasnosels'ki\u{\i}-Mann theorem \cite{Bauschke}.

\subsection{Non-Expansivity of $U_{f,g}$ and $U_{g,f}$}
In this subsection, we obtain a condition which ensures that $U_{f,g}$ and $U_{g,f}$ are non-expansive.
\begin{defn}
An operator $S:\mathbb{R}^n\to \mathbb{R}^n$ is said to be non-expansive if
\begin{equation}
\|S(x) - S(y)\|_2 \leq \|x - y\|_2,
\end{equation}
for all $x$, $y$.
\end{defn}
For this purpose, we separately study the operators
\begin{subequations}
\begin{align}
U_{f} &= 2J_{\alpha f} - I, \\
U_{g} &= 2J_{\alpha g} - I.
\end{align}
\end{subequations}
Note that using $U_f$ and $U_g$, the operators in \eqref{eqn:Ufgf} can be expressed as $U_{f,g} = U_f\,U_g$ and $U_{g,f} = U_g\,U_f$. We start our study with $U_f$.

\subsubsection{A Lipschitz Constant for $U_f$} 
\begin{lem}\label{lem:f}
Suppose $f$ is a second-order differentiable function such that $f(x) - \frac{\rho}{2}\|x\|_2^2$ is convex and $\nabla f$ is Frechet differentiable with Lipschitz constant $\sigma$, i.e., 
\begin{equation}
\| \nabla f(x) - \nabla f(y) \|_2 \leq \sigma \|x - y \|_2.
\end{equation}
Then, ${U_f = 2J_{\alpha f} - I}$\, satisfies
\begin{equation}\label{eqn:concL1}
\|U_f(x) - U_f(y)\|_2 \leq \max \left( \frac{| 1 - \alpha\, \sigma |}{1 + \alpha\, \sigma}, \frac{|1 - \alpha \,\rho | }{1 + \alpha\,\rho} \right)\,\|x-y\|_2.
\end{equation}

\begin{proof}
Let $F = \nabla f$. Since $F$ is single-valued and $f$ is convex, it follows that $G_{\alpha} = I + \alpha \,F$ is a bijection on $\R^n$ by Minty's theorem (see Thm.~21.1 in \cite{Bauschke}). Now, let $d_{\alpha}$ denote the Frechet derivative of $G_{\alpha}$. It follows by the assumptions on $f$ that the spectrum of $d_{\alpha}$ at any $x\in \R^n$ is contained in $[1 + \alpha\,\rho,1+\alpha\,\sigma]$. It also follows from the inverse function theorem (see e.g. Thm. 16.12 in \cite{Fitzpatrick}) that $J_{\alpha\,f} = G_{\alpha}^{-1}$ is continuous, differentiable and the spectrum of the Frechet derivative of $J_{\alpha f}$ is contained in $[(1+\alpha\,\sigma)^{-1}, (1+\alpha\,\rho)^{-1}]$. As a consequence, the spectrum of the Frechet derivative of $U_f$ is contained in 
\begin{equation}
\left[ \frac{1 - \alpha\,\sigma}{1+\alpha\,\sigma},  \frac{1 - \alpha\,\rho}{1+\alpha\,\rho} \right],
\end{equation}
from which, \eqref{eqn:concL1} follows.
\end{proof}
\end{lem}

\subsubsection{A Lipschitz Constant for $U_g$}

We will make use of the fundamental result below, which can be found in \cite{Bauschke}.
\begin{lem}\label{lem:fn}
Suppose $q$ is a convex function. Then, $J_{\alpha q}$ can be written as,
\begin{equation}
J_{\alpha q} = \frac{1}{2} (I +  S),
\end{equation}
where $S$ is non-expansive.
\end{lem}

\begin{lem}\label{lem:g}
Suppose $g$ is a $\rho$ weakly convex function and ${\alpha < 1/\rho}$. Also, let ${U_g = 2J_{\alpha g} - I}$. Then,
\begin{equation}
\|U_g(x) - U_g(y)\|_2 \leq \frac{1 + \alpha\,\rho}{ 1- \alpha\,\rho} \, \|x - y \|_2.
\end{equation}
\begin{proof}
Since $g$ is $\rho$-weakly convex, the function ${\tilde{g}(x) = g(x) + \frac{\rho}{2} \|x\|_2^2}$ is convex and $J_{\gamma \tilde{g}} = \frac{1}{2} (I + S)$ for a non-expansive $S$, by Lemma~\ref{lem:fn}. Now note that,
\begin{align}
J_{\gamma \tilde{g}}(x) &= \arg \min_z \frac{1}{2\gamma} \|z - x \|_2^2 + \frac{\rho}{2} \|z\|_2^2+ g(z)  \\
 &= \arg \min_z \frac{1}{2} \left\|z - \frac{1}{1 + \gamma \rho} x \right\|_2^2 + \frac{\gamma}{1+\gamma\rho}g(z)\\
&= J_{\gamma(1 + \gamma\rho)^{-1}g} \left( \frac{1}{1 + \gamma \rho} x \right)
\end{align}
Substituting $\alpha = \gamma/(1 + \gamma \rho)$, we can thus write
\begin{align}
J_{\alpha g} \left(  x \right) &= J_{\gamma \tilde{g}} \left( \frac{1}{1 - \alpha \rho} x \right)\\
&=\frac{1}{2}  \frac{1}{1 - \alpha \rho} x + \frac{1}{2} S\left( \frac{1}{1 - \alpha \rho} x \right).
\end{align}
We now have
\begin{align}
\| U_g (x) &- U_g (y) \|_2 
= \left\| \frac{\alpha\,\rho}{1 - \alpha\,\rho} (x -y )\,  \right. \nonumber \\
& \quad \left.+ S\left( \frac{1}{1 - \alpha \rho} x \right) -  S\left( \frac{1}{1 - \alpha \rho} y \right)\right\|_2\\
&\leq \frac{\alpha\,\rho}{1 - \alpha\,\rho} \| x -y \|_2 \nonumber 
\\& \quad + \left\| S\left( \frac{1}{1 - \alpha \rho} x \right) -  S\left( \frac{1}{1 - \alpha \rho} y \right)\right\|_2\\
&\leq \frac{1 + \alpha\,\rho}{1 - \alpha\,\rho} \| x -y \|_2.
\end{align}
\end{proof}
\end{lem}

We are now ready for the main result of this subsection.
\begin{prop}\label{prop:NonExpansive}
If 
\begin{equation}
\alpha \leq \frac{1}{\sqrt{\rho\,\sigma}},
\end{equation}
then $U_{fg}$ and $U_{gf}$ are non-expansive.
\begin{proof}
Combining Lemmas~\ref{lem:f}, \ref{lem:g}, for $\alpha < 1/\rho$, we see that $U_{fg}$ and $U_{gf}$ are non-expansive if
\begin{equation}
\max \left( \frac{| 1 - \alpha\, \sigma |}{1 + \alpha\, \sigma}, \frac{1 - \alpha \,\rho }{1 + \alpha\,\rho} \right)\,\times \frac{1 + \alpha\,\rho}{1 - \alpha\,\rho} \leq 1.
\end{equation}
This is equivalent to
\begin{equation}
\max \left( \frac{| 1 - \alpha\, \sigma |}{1 + \alpha\, \sigma}, \frac{1 - \alpha \,\rho }{1 + \alpha\,\rho} \right)\,\leq  \frac{1 - \alpha\,\rho}{1 + \alpha\,\rho}.
\end{equation}
This in turn is satisfied if
\begin{equation}
\frac{| 1 - \alpha\, \sigma |}{1 + \alpha\, \sigma} \leq \frac{1 - \alpha \,\rho }{1 + \alpha\,\rho}.
\end{equation}
Noting that $\sigma \geq \rho$, this is found to be equivalent to
\begin{equation}
\alpha \leq \frac{1}{\sqrt{\rho\,\sigma}}.
\end{equation}
\end{proof}
\end{prop}

\subsection{Fixed Points of $U_{f,g}$ and $U_{g,f}$}
In this subsection, we study the fixed points of $U_{f,g}$ and $U_{g,f}$ and establish their relation with the minimizers of $h$.

Let us first recall a result from \cite{bay15arXiv}.
\begin{prop}(Prop.~10 in \cite{bay15arXiv})\label{prop:fixed1}
Suppose the conditions stated in Prop.~\ref{prop:main} hold. Then,
$x^*$ minimizes $f+g$ if and only if ${x^* = J_{\alpha\,g}(I - \alpha\,\nabla f) (x^*)}$.
\end{prop}

\begin{prop}\label{prop:fixed2}
Suppose the conditions stated in Prop.~\ref{prop:main} hold.
\begin{enumerate}[(a)]
\item $z = U_{f,g}(z)$  if and only if $J_{\alpha g}(z)$ minimizes $f+g$.
\label{Ufg}
\item $z = U_{g,f}(z)$  if and only if $J_{\alpha f}(z)$ minimizes $f+g$. \label{Ugf}
\end{enumerate}
\begin{proof}
Using Prop.~\ref{prop:fixed1}, we find that $x^*$ minimizes $f+g$ if and only if
\begin{equation}
(I + \alpha\,\nabla f)(x^*)  = \bigl(J_{\alpha\,g}(I - \alpha\,\nabla f) + \alpha\,\nabla f \bigr)(x^*).
\end{equation}
Now let $z = (I + \alpha\,\nabla f)(x^*)$ so that  $x^* =J_{\alpha f}(z)$. Then, $x^*$ minimizes $f+g$ if and only if
\begin{equation}
z  = \bigl(J_{\alpha\,g}(2 J_{\alpha\,f} - I)  + (I - J_{\alpha\, f}) \bigr)(z).
\end{equation}
But this is equivalent to
\begin{equation}
z  = (2 J_{\alpha\,g} - I)\,(2 J_{\alpha\,f} - I)(z).
\end{equation}
Thus follows \eqref{Ugf}.

Now suppose 
\begin{equation}\label{eqn:z}
z  = (2 J_{\alpha\,g} - I)\,(2 J_{\alpha\,f} - I)(z).
\end{equation}
and $x^* = J_{\alpha f}(z)$ or $z = (I + \alpha\,\nabla f)(x^*)$.
Let us define
\begin{equation}
q = (2 J_{\alpha\,f} - I)\,(z) = 2x^* - z.
\end{equation}
We have, by \eqref{eqn:z},
\begin{equation}
q  = (2 J_{\alpha\,f} - I)\,(2 J_{\alpha\,g} - I)\,(q).
\end{equation}
and $z = (2 J_{\alpha\,g} - I)\,(q)$. Thus,
\begin{equation}
x^* = \frac{q+z}{2} = \frac{q + (2 J_{\alpha\,g} - I)(q)}{2} = J_{\alpha\,g}\,(q).
\end{equation}
Thus follows \eqref{Ufg}.
\end{proof}
\end{prop}

\subsection{Proof of Prop.~\ref{prop:main}}
The final ingredient for the proof of Prop.~\ref{prop:main} is the Krasnosels'ki\u{\i}-Mann Theorem.
\begin{prop}\label{prop:KM} \cite{Bauschke}
Suppose $U$ is non-expansive, its set of fixed points is non-empty and ${S = (1-\lambda)\,I+ \lambda\,U}$ for $0< \lambda< 1$. Also, let the sequence $z^n$ be defined as $z^{n+1} = S\,(z^n)$. Then, the sequence $z^n$ converges to a fixed point of $U$. \qed
\end{prop}

We are now ready for the proof of Prop.~\ref{prop:main}. Observe that the algorithm in \eqref{eqn:DR1a} sets $z^{n+1} = S(z^n)$, where ${S = (1-\lambda)\,I + \lambda\,U_{f,g}}$ with $0< \lambda< 1$. But $U_{f,g}$ was shown to be non-expansive in Prop.~\ref{prop:NonExpansive}. Therefore, by Prop.~\ref{prop:KM}, $z^n$'s converge to a fixed point of $U_{f,g}$ (which exist by assumption), say $z^*$. Prop.~\ref{prop:fixed2}~\ref{Ufg} notes that if $z^*$ is a fixed point of $U_{f,g}$, then $J_{\alpha\,g}(z^*)$ minimizes $f + g$, which proves the claim in Prop.~\ref{prop:main}\ref{alg:a}.

For part \eqref{alg:b} of Prop.~\ref{prop:main}, the argument is similar. We only need to replace $U_{f,g}$ with $U_{g,f}$ and refer to Prop.~\ref{prop:fixed2}\ref{Ugf} instead of Prop.~\ref{prop:fixed2}\ref{Ufg}.

\subsection{Discussion : Splittings Allowed by Prop.~\ref{prop:main}}
The algorithms in Prop.~\ref{prop:main} do not allow $f$ to be non-smooth. This can be constraining in some scenarios. For instance, for  $K$ a subspace of $\R^n$, consider the problem,
\begin{equation}\label{eqn:prob}
\min_{x\in \R^n} \frac{1}{2} \|y - x\|_2^2 + i_K(x) + \phi(x),
\end{equation}
where $i_K(x)$ is the characteristic function of $K$ and $\phi$ is 1-weakly convex. Note that the cost function in \eqref{eqn:prob} is convex. Now if we set $f(x) = \frac{1}{2} \|y - x \|_2^2 + i_K(x)$, then  ${ J_{\alpha\,f}(z) = P_K(z + \alpha\,y)}$, where $P_K$\, denotes the (linear) projection operator onto $K$. It follows that
\begin{equation}
(2 J_{\alpha f} - I )x_0 - (2 J_{\alpha f} - I )x_1 = (2P_K - I)\,(x_0 - x_1).
\end{equation}
But the operator $(2P_K - I)$ is not a contraction, for any value of $\alpha$. Therefore, $U_{f,g} = (2 J_{\alpha f} - I )\,(2 J_{\alpha g} - I )$ may not be non-expansive in this case. But this means that we cannot employ Prop.~\ref{prop:KM}. Consequently, the proof of Prop.~\ref{prop:main} outlined above does not cover this case.

We remark that, in conrast to the algorithms in Prop.~\ref{prop:main}, the algorithms in Prop.~\ref{prop:shift}  would employ a splitting as $f(x) = i_K(x)$, $g(x) = \frac{1}{2} \|y -x \|_2^2 + \phi(x)$ and avoid working with the proximal mapping associated with the non-convex $\phi$.

\section{Shifting the Quadratic (Proof of Prop.~\ref{prop:shift})}\label{sec:shift}
In this section, we present another approach for handling a weakly convex penalty function. Prop.~\ref{prop:shift} follows as a consequence of this discussion.

We start by rewriting the minimization problem in \eqref{eqn:problem} as 
\begin{equation}
\min_x \underbrace{\left( f(x) - \frac{\rho}{2}\, \|x\|_2^2 \right)}_{\tilde{f}}+ \underbrace{\left( g(x) + \frac{\rho}{2}\, \|x\|_2^2 \right)}_{\tilde{g}}.
\end{equation}
Thanks to the properties of $f$ and $g$, it follows that $\tilde{f}$ and $\tilde{g}$ are convex. Thus we can apply the regular Douglas-Rachford algorithm for this new splitting. This leads to iterations of the form
\begin{equation}\label{eqn:DR2}
z^{n+1} = \Bigl( (1-\lambda) I + \lambda\,\bigl(2J_{\alpha \tilde{f}} - I \bigr) \, \bigl(2J_{\alpha \tilde{g}} - I \bigr) \Bigr)\,(z^n).
\end{equation}
Let us now find a relation between $J_{\alpha \tilde{f}}$, $J_{\alpha \tilde{g}}$ and the proximity operators of $f$ and $g$. Observe that 
\begin{subequations}
\begin{align}
J_{\alpha \tilde{g} }(x) &= \arg \min_z \frac{1}{2 \alpha}\|z -x \|_2^2 + \frac{\rho}{2} \|z\|_2^2 +  g(z) \\
&=\arg \min_z \frac{1}{2} \left \|z - \frac{1}{1+ \alpha\,\rho}x \right \|_2^2 + \frac{\alpha}{1+ \alpha\,\rho}\, g(z) \\
&= J_{\beta_1 g} \left( \frac{\beta_1}{\alpha} x \right), \label{eqn:K1alt}
\end{align}
\end{subequations}
where $\beta_1 = \alpha/ (1+\alpha\,\rho)$.

Now suppose that $\alpha\,\rho < 1$. We similarly obtain
\begin{subequations}
\begin{align*}
J_{\alpha \tilde{f} }(x) &= \arg \min_z \frac{1}{2 \alpha}\|z -x \|_2^2 - \frac{\rho}{2} \|z\|_2^2 +  f(z) \\
&=\arg \min_z \frac{1}{2} \left \|z - \frac{1}{1 -  \alpha\,\rho}x \right \|_2^2 + \frac{\alpha}{1- \alpha\,\rho}\, f(z) \\
&= J_{\beta_2 f} \left( \frac{\beta_2}{\alpha} x \right), \label{eqn:K2alt}
\end{align*}
\end{subequations}
where $\beta_2 = \alpha/ (1 - \alpha\,\rho)$.

Plugging these in \eqref{eqn:DR2}, we obtain the first algorithm in Prop.~\ref{prop:shift}. The second algorithm is obtained by swapping  $J_{\alpha \tilde{f}}$ and $J_{\alpha \tilde{g}}$. Since $\tilde{f}$ and $\tilde{g}$ are convex, the convergence claims in Prop.~\ref{prop:shift} follow as a consequence of Prop.~\ref{prop:convex}.

\section{Convergence Speed of the Algorithms}\label{sec:speed}

We provided two sets of algorithms in Prop.~\ref{prop:main} and Prop.~\ref{prop:shift} and noted in the Introduction that the algorithms in Prop.~\ref{prop:shift} apply to a wider class of splittings. However, the algorithms in Prop.~\ref{prop:main} are still of interest because they can be faster under certain conditions. In this section, we provide some theoretical evidence to support this claim. Specifically, we study the contraction rates of operators associated with Peaceman-Rachford  iterations (i.e., $\lambda = 0$ in \eqref{eqn:DR})  for the two algorithms (see $T_\alpha$ and $V_{\alpha}$ in \eqref{eqn:UV}), under a strong convexity assumption. Althuogh these contraction rates do not automatically transfer to Douglas-Rachford iterations, they do provide some intuition on how the convergence speed may vary with the properties of the problem. This intuition is also supported by numerical experiments in Sec.~\ref{sec:exp}.

For this section, we make the following assumptions.
\begin{enumerate}[$(\text{A}_1)$]
\item $g$ is $\rho$-weakly convex,
\item $f$ is second order differentiable and $\nabla f$ satisfies ${\| \nabla f(x) - \nabla f(y) \|_2 \leq \sigma \|x - y\|_2}$,
\item for some $s$ with $s > \rho$, the function $f(x) - \dfrac{s}{2}\|x\|_2^2$ is convex.
\end{enumerate}
Under these assumptions, we will see that the operators $T_{\alpha}$, $V_{\alpha}$ extracted from \eqref{eqn:DR1a}, \eqref{eqn:DRKa} as
\begin{subequations}\label{eqn:UV}
\begin{align}
T_{\alpha} &= \bigl(2J_{\alpha f} - I \bigr) \, \bigl(2J_{\alpha g} - I \bigr),\\
V_{\alpha} &= \bigl(2K_2 - I \bigr) \, \bigl(2 K_1 - I \bigr),
\end{align}
\end{subequations}
are not only non-expansive but they are actually contractions \cite{Bauschke}. We can therefore compare their contraction rates to gain some idea about their convergence speed. 

\begin{prop}\label{prop:speed}
Suppose the assumptions $A_1$, $A_2$, $A_3$ are in effect. Also, let $T_\alpha$, $V_{\alpha}$ be defined as in \eqref{eqn:UV}.

If $\alpha \leq 1/\sqrt{\sigma\,s}$, then
\begin{equation}\label{eqn:c1}
\| T_{\alpha}(x) - T_{\alpha}(y) \|_2 \leq \frac{(1 - \alpha^2 s\,\rho) - \alpha(s-\rho)}{(1 - \alpha^2 s\,\rho) + \alpha(s-\rho)}  \|x - y \|_2.
\end{equation}

If $\alpha \leq 1/s$, then
\begin{equation}\label{eqn:c2}
\| V_{\alpha}(x) - V_{\alpha}(y) \|_2 \leq  \gamma  \| x - y \|_2,
\end{equation}
where 
\begin{equation}
\gamma = \max\left( \frac{|1 - \alpha(\sigma-\rho)|}{1 + \alpha(\sigma-\rho)}, \frac{1  - \alpha(s-\rho)}{1 + \alpha(s-\rho)} \right).
\end{equation}

\begin{proof}
Suppose $\alpha \leq 1/\sqrt{\sigma\,s}$.
This implies,
\begin{equation}
\frac{\alpha\, \sigma - 1}{1 + \alpha\, \sigma} \leq  \frac{1 - \alpha \,s  }{1 + \alpha\,s}.
\end{equation}
Using Lemma~\ref{lem:f}, we obtain,
\begin{equation}
\| (2J_{\alpha\,f} - I)x - (2J_{\alpha\,f} - I)y \|_2 \leq \frac{1 - \alpha\,s}{1+ \alpha \,s} \| x - y\|_2.
\end{equation}
Since $s > \rho$, we also have $\alpha \leq 1/\sqrt{\sigma\,\rho}$.  From Lemma~\ref{lem:g}, we now obtain, 
\begin{equation}
\| (2J_{\alpha\,g} - I)x - (2J_{\alpha\,g} - I)y \|_2 \leq \frac{1  + \alpha\,\rho}{1 -  \alpha \,\rho} \| x - y\|_2.
\end{equation}
Combining the two, we obtain
\begin{align}
\| T_{\alpha}(x) - T_{\alpha}(y) \|_2 & \leq  \frac{1 - \alpha\,s}{1+ \alpha \,s}\, \frac{1  + \alpha\,\rho}{1 -  \alpha \,\rho}\, \| x - y\|_2\\
&=\frac{(1 - \alpha^2 s\,\rho) - \alpha(s-\rho)}{(1 - \alpha^2 s\,\rho) + \alpha\,(s-\rho)}  \|x - y \|_2.
\end{align}

Let us now show the second part. Assume that $\alpha \leq 1/s$. By the discussion in the proof of Lemma~\ref{lem:f}, the spectrum of the Frechet derivative of $J_{\beta_2f}(x)$ lies in 
\begin{equation}
\left[ \frac{1 - \beta_2\,\sigma}{1+\beta_2\,\sigma},  \frac{1 - \beta_2\,\rho}{1+\beta_2\,\rho} \right].
\end{equation}

Since $K_2(x) = J_{\beta_2 f}\,\left( \dfrac{\beta_2}{\alpha} x \right)$, and $\beta_2 = \alpha/(1 - \alpha\,\rho)$,
it follows that the spectrum of the Frechet derivative of $2K_2 - I$ lies in 
\begin{equation}
\left[ \frac{1 - \alpha (\sigma - \rho)}{1 + \alpha (\sigma - \rho) }, \frac{1 - \alpha (s - \rho)}{1 + \alpha (s - \rho )} \right].
\end{equation}
Also, since $K_1 = J_{\alpha \tilde{g} }$, where $\tilde{g}(x) = g(x) + \dfrac{\rho}{2}\|x\|_2^2$ is convex, it follows by  Lemma~\ref{lem:fn} that ${2 K_1 - I}$ is non-expansive. Combining these observations, \eqref{eqn:c2} follows.
\end{proof}
\end{prop}

We remark that if the same choice of the step parameter $\alpha$ is used, in the worst case, repeated applications of $T_{\alpha}$ is likely to converge faster than $V_{\alpha}$ since 
\begin{equation}
\frac{(1 - \alpha^2 s\,\rho) - \alpha(s-\rho)}{(1 - \alpha^2 s\,\rho) + \alpha(s-\rho)}  \leq  \frac{1  - \alpha(s-\rho)}{1 + \alpha(s-\rho)}.
\end{equation}
However, the contraction rates are inversely proportional to $\alpha$ and the highest value of $\alpha$ allowed by the two algorithms are different in Prop.~\ref{prop:speed}. Setting $\alpha = 1/\sqrt{\sigma\,s}$, we find the minimum contraction rate of $T_{\alpha}$ as 
\begin{equation}
\mu_T = \frac{(1 - \eta)\sqrt{\gamma} - (\gamma - \eta)}{(1-\eta)\sqrt{\gamma} +  (\gamma - \eta)},
\end{equation}
where $\gamma = s/\sigma$ and $\eta = \rho / \sigma$.
Observe that as $s \to \sigma$, that is, as $\gamma \to 1$, $\mu_T$ decreases to zero. 
Therefore the convergence speed in this case depends on the problem.
On the other hand, for $\alpha = 1/s$, we find that 
\begin{equation}
\frac{\eta}{2-\eta} \leq \frac{1  - \alpha(s-\rho)}{1 + \alpha(s-\rho)}.
\end{equation}
Therefore, the contraction rate of $V_{\alpha}$ is bounded below by $\eta/(2-\eta)$. Since this expression is independent of $\gamma$, the convergence speed of repeated applications of $V_{\alpha}$ is not affected by the conditioning of the smooth term.

To summarize, based on this discussion, we expect the algorithm in Prop.~\ref{prop:main} to converge faster than the algorithm in Prop.~\ref{prop:shift} for problems where $\nabla^2 f$ is better approximated by $I$, so that the ratio $s/\sigma$ is close to unity. We demonstrate this with a numerical experiment in Exp.~\ref{exp:2} below.

\section{Experiments}\label{sec:exp}

In this section, we demonstrate that the proposed algorithms do converge and provide an empirical comparison of their convergence rate. Matlab code for the experiments can be found at \texttt{http://web.itu.edu.tr/ibayram/NCDR/}.

We experiment with a simple sparse signal recovery problem. We use $90\times 120$ convolution matrices $H$  associated with invertible filters (so that $H^T\,H$ is invertible) to construct the observed signal. Using a sparse $x$, we produce the observed signals as 
\begin{equation}
y = H\,x + u,
\end{equation}
where $u$ denotes white Gaussian noise. We set the signal variance so that the SNR is 10dB. 

\begin{figure}
\centering
\includegraphics[scale=1]{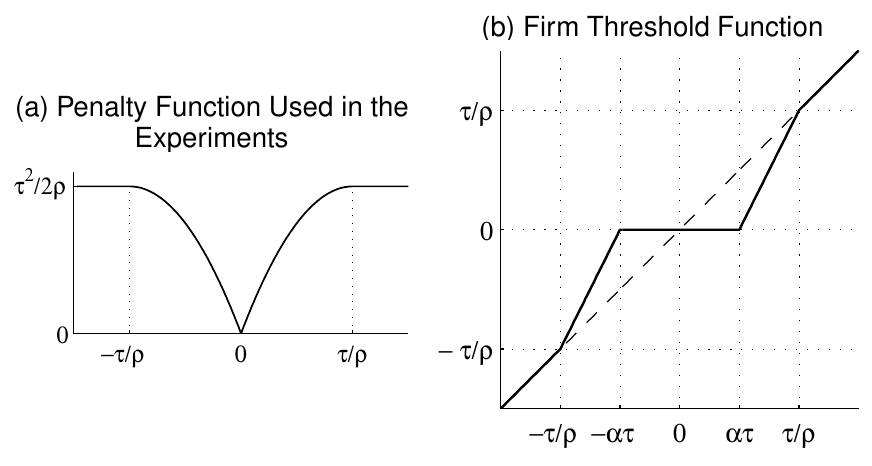}\\
\caption{The penalty and threshold function used in the experiments. \label{fig:DemoTHold}}
\end{figure}

In both of the experiments below, we use a separable penalty function based on $P:\mathbb{R} \rightarrow \mathbb{R}$ given as (see Fig.~\ref{fig:DemoTHold}a)
\begin{equation}
P_{\tau,\rho}(t) = \begin{cases}
\tau\,|t| - t^2/(2\rho), &\text{ if } |t| < \tau/\rho,\\
\tau^2/(2\rho), &\text{ if } |t| \geq \tau/\rho,\\
\end{cases}
\end{equation}
This  function is $\rho$-weakly convex and is the penalty function associated with `firm-thresholding' \cite{gao97p855}. The proximity operator for $P$, namely the firm-threshold, is given by (provided $\alpha\,\rho < 1$),
\begin{equation}\label{eqn:DemoTHold}
T_{\alpha}(t) = \begin{cases}
0, &\text{if } |t| < \alpha\,\tau,\\
(1-\alpha\,\rho)^{-1}\,(t - \alpha\,\tau), &\text{if } \alpha\,\tau \leq |t| < \tau/\rho,\\
t, &\text{if } \tau/\rho \leq |t|.
\end{cases}
\end{equation}
$T_{\alpha}(t)$ is depicted in Fig.~\ref{fig:DemoTHold}b. Observe that for ${\alpha \,\tau < |s| < \tau/\rho}$, the derivative of the threshold function exceeds unity. Therefore the threshold function is not non-expansive.

In the setup described above, we let $s$ and $\sigma$ denote the least and the greatest eigenvalue of $H^T\,H$. We also set ${\tau = 3\,\rho \,\std(u)}$, where $\std(u)$ denotes the standard deviation of noise. We obtain the estimate of the sparse signal as,
\begin{equation}\label{eqn:ExpCost}
x^* = \arg \min_t \frac{1}{2} \|y - H\,t\|_2^2 + \sum_i P_{\tau,\rho}(t_i).
\end{equation}

\begin{experiment}\label{exp:1}

In the first experiment, we take $\rho = s$ and pick a filter such that 
the ratio $\sigma / s$ is 15.96.

\begin{figure}
\centering
\includegraphics[scale = 1]{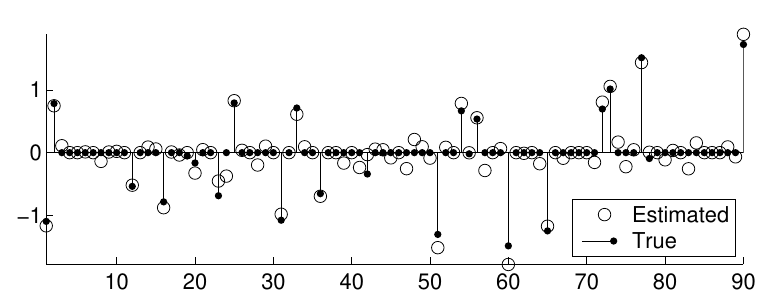}\\
\caption{The underlying unknown sparse signal and its estimate from Experiment~\ref{exp:1}.}\label{fig:signals}
\end{figure}

\begin{figure}
\centering
\includegraphics[scale = 1]{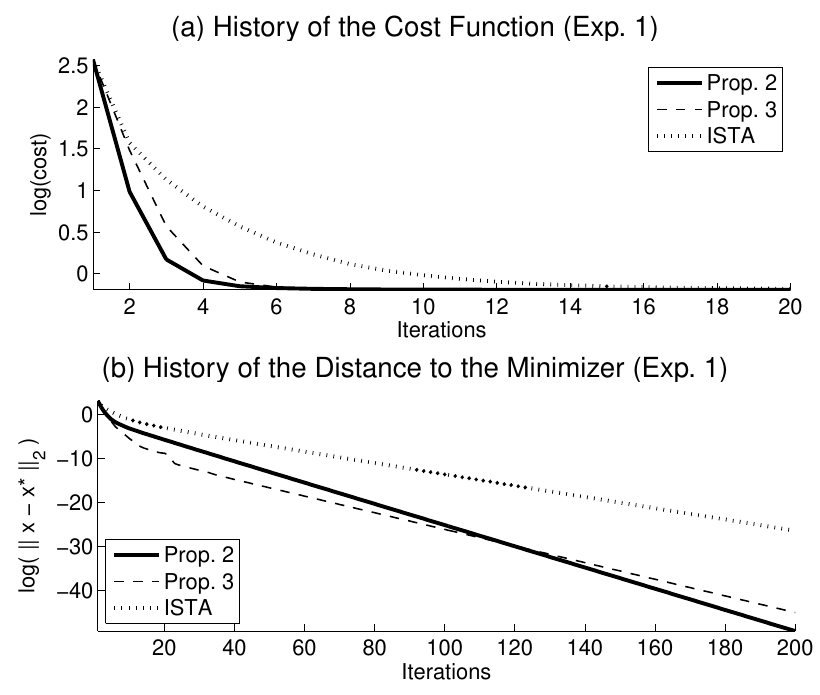}\\
\caption{History of the (a) cost and  (b) distance to the minimizer for the algorithms from Experiment~\ref{exp:1}}\label{fig:Cost}
\end{figure}

The history of the cost function with iterations is shown in Fig.~\ref{fig:Cost}a for the two algorithms from Propositions~\ref{prop:main} and \ref{prop:shift}. For both algorithms, we used a value $\alpha$ that is close to the allowed upper bound. 

In order to assess the convergence speed, we first obtain an estimate of the minimizer $x^*$ and then monitor the distance to this estimate with iterations. In order to prevent any bias, we obtain the estimate of the minimizer by another algorithm, namely the iterative shrinkage/thresholding algorithm (ISTA) for 10K iterations (see \cite{bay15arXiv} for a discussion of ISTA for weakly convex penalties). The logarithm of the Euclidean distance to $x^*$ with respect to iterations is shown in Fig.~\ref{fig:Cost}b for the two algorithms in Prop.~\ref{prop:main} and Prop.~\ref{prop:shift}. We have observed that both versions of the DR algorithm perform faster than ISTA consistently, in terms of reducing the cost as well as approaching the minimizer. Depending on the realization of the random noise and the underlying sparse signal, the convergence speed of the algorithms in Prop.~\ref{prop:main} and Prop.~\ref{prop:shift} vary. Although both algorithms converge faster than ISTA, we did not observe a meaningful trend so as to conclude that one of them converges faster. 

\end{experiment}

\begin{experiment}\label{exp:2}
In this experiment, we demonstrate the discussion in Section~\ref{sec:speed}. For that, we chose the filter associated with $H$ so that the ratio of the greatest and smallest eigenvalues of $H^T\,H$ is  $\sigma / s = 5.44$. 
We also set $\rho = s / 2$. As in Experiment~\ref{exp:1}, we set $\alpha$ to be close to the allowed upper bound. Note that for this problem, the Hessian of the data fidelity term, namely $\nabla^2 f(x)$, is closer to the identity than that in Experiment~\ref{exp:1}. Therefore, we expect the first algorithm to perform faster, by the discussion at the end of in Section~\ref{sec:speed}.

The logarithm of the Euclidean distance to the minimizer with respect to iterations is shown in Fig.~\ref{fig:dist2}. We observed consistently in this case that the first algorithm has a higher convergence rate, in agreement with our expectation.

\begin{figure}
\centering
\includegraphics[scale = 1]{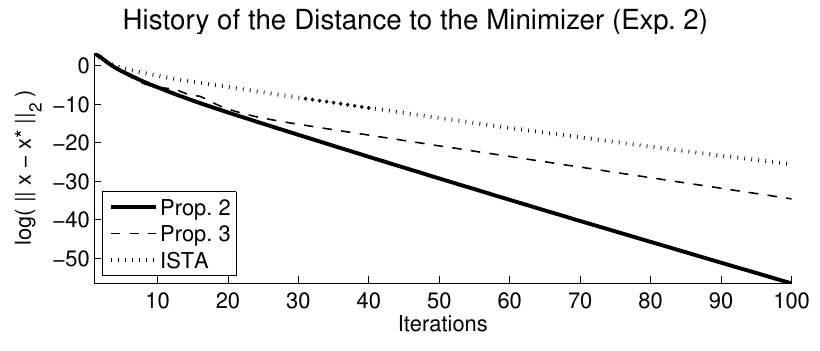}
\caption{History of the distance to the minimizer for the algorithms from Experiment~\ref{exp:2}}\label{fig:dist2}
\end{figure}

\end{experiment}

\section{Outlook}\label{sec:conc}
Convex minimization problems that contain non-convex terms are receiving more attention in sparse signal processing due to their attractive property of preserving convexity while enforcing sparsity stronger than convex alternatives like the $\ell_1$ norm or mixed norms. In this paper, we considered the convergence of the Douglas-Rachford algorithm in such a setting. We showed that the original Douglas-Rachford iterations do converge even though some proximity operators are associated with non-convex (but weakly-convex) functions, provided that a condition on the `step-size' is satisfied.

A possible extension of the current work is to consider a scenario where the sum to be minimized contains  more than two functions. In the convex case, such an extension can be handled by the parallel proximal algorithm (PPXA)  \cite{combettes_chp} derived from the Douglas-Rachford algorithm. However, in a preliminary study, we observed that the analysis in this paper does not easily generalize to show that PPXA converges when some of the terms are weakly-convex. Although an algorithm can be derived by `shifting the quadratic' as we did in Prop.~\ref{prop:shift}, it is of interest to find conditions under which algorithms applicable in the convex case can also be used in this setting, without modification.  We hope to study this problem in the near future.

\bibliographystyle{plain}

\end{document}